\documentclass[12pt,leqno]{article}
\usepackage{amsmath}
\usepackage{amsfonts}
\newcounter{ss}
\newcounter{sss}          

\newcommand{\RR}{\Bbb R}

\newcommand{\CC}{\Bbb C}
\newcommand{\ZZ}{\Bbb Z}
\newcommand{\NN}{\Bbb N}

\newtheorem{theorem}{Theorem}[section]

\newtheorem{definition}{Definition}[section]

\newtheorem{lemma}{Lemma}[section]
\newtheorem{corollary}{Corollary}[section]

\begin{document}

\begin{center}
{\bf\LARGE Equifocal families in symmetric spaces \\
  of compact type}
\end{center}
\begin{center}
{\Large Martina Br\"uck}
\end{center}
 
\begin{center}
{\bf Abstract}
\end{center}
\vspace{-0.6cm} 
\begin{quote}
\small An equifocal submanifold $M$ of a symmetric space $N$ of compact type induces a  foliation
with singular leaves on $N$. In this paper we will show how to reconstruct the equifocal foliation starting from
one of  the singular leaves, the so-called {\it focal
manifolds}. To be more concrete: The equifocal submanifold is equal to a
partial tube $\tilde{B}_{\xi}$ around the focal manifold and we
will show how to construct $\tilde{B}_{\xi}$ in this paper. Moreover, we will find a geometrical
characterization of focal manifolds.
\end{quote}

\section{Introduction}

   The classification of group actions on a given manifold and
the description of the geometry of their orbits are important problems in Riemannian geometry. As the set
of different  group
actions is huge, one usually has to restrict the considerations to  certain classes of group actions,    for example by
assuming some ``nice" properties of the orbits. In this paper we want to study  {\it hyperpolar actions} on symmetric
spaces    and   {\it equifocal}  (singular)  foliations. 
 An operation on a symmetric space is called {\it hyperpolar } if it admits
a flat section. This means,
 there exists a  flat, totally geodesic, closed subspace of the ambient
space that meets every orbit and meets it orthogonally at every intersection point.   A certain class of hyperpolar actions
was introduced by   Hermann   in \cite{Her}: Let
$(G,K)$ and $(G,\tilde{K})$ be two Riemannian symmetric pairs of compact type. Then, $\tilde{K}$ acts hyperpolarly on 
$G/K$.\renewcommand{\thefootnote}{}

  Let\footnote{1991 {\it Mathematics Subject Classification.} Primary 53C40; Secondary 53C35.}
 $M\subset N$ be a connected Riemannian submanifold of a symmetric space 
and  $\exp^{\bot}:=(\exp^N |\bot M)$   the restriction of 
$\exp^N$ to the {\it normal bundle} $\bot M$ of $M\subset N$. 
A vector $v \in \bot_p M$ is called a  {\it  focal normal of multiplicity }
  $l>0$ if the kernel of   $\exp^{\bot} _{\ast |v } $ is $l$-dimensional. 
$\exp ^{\bot} (v)$ is then called a {\it
focal point of multiplicity $l$}.
We say that  $M$ has {\it constant focal distances} if for every curve
$c:[0,1]\rightarrow M$ and every   parallel normal vector field $
w$ along\footnote{Supported by     Deutscher Akademischer Austauschdienst (DAAD)}
$c$ the kernel of the differential of $\exp^{\bot}$ is constant along
$w$, that is to say  dim(ker( $\exp ^{\bot }_{\ast | w(t)}))= $constant. 
  A submanifold is called {\it equifocal }  if it has  {\it abelian}
 and {\it globally flat normal bundle }  and
{\it constant focal distances.} Here, we say that  $M$ has {\it  abelian normal bundle }   if
$\exp^N (\bot_p M)$ is contained in a flat, totally geodesic submanifold for all $p\in M$. We say that  $M$ has 
{\it globally flat normal bundle }  if the   normal holonomy group $\phi_p$ at any point $p\in M$ is trivial. And we say
that
$M$ has {\it  flat normal bundle } if the   reduced normal holonomy group  $\phi_p^{\ast}$ is trivial. 

  Terng and Thorbergsson showed in \cite{TT} (compare Theorem 2.1) that the principal orbits of  hyperpolar 
actions on symmetric spaces of compact type are   equifocal.   If $\tilde{M}\subset N$ is an embedded equifocal submanifold
of a symmetric space of compact type and $\tilde{\xi} :\tilde{M}\rightarrow\bot \tilde{M}$ is a globally defined  
parallel  normal vector field then
$\exp^N(\tilde{\xi} (\tilde{M}))=:M$ is an embedded submanifold of $N$
 (compare
\cite{TT}). If $\mbox{dim}(\tilde{M})=\mbox{dim}(M)$ then  $M$
is {\it locally equifocal} (maybe the normal
holonomy group is not trivial, but still it is discrete). If
$\mbox{dim}(\tilde{M})>  \mbox{dim}(M)$ then $M$ is called a {\it focal manifold of  $\tilde{M}$}.  The singular orbits of 
hyperpolar actions on symmetric spaces of compact type are focal manifolds of the principal ones. Up to now
only little was
known about the properties of the singular orbits or more generally: about the focal manifolds of equifocal ones. And it
was not clear if and in which way the equifocal manifolds  can be reconstructed starting from a given focal manifold.
In this paper we will prove:

    The equifocal submanifold is equal to a partial tube
$\exp^N(B_{\xi})=\tilde {B}_{\xi}$ around the focal manifold. The fibres of $B_{\xi}$  are
homogeneous. To be more concrete  they are orbits of a Lie subgroup $\hat{G}_p \subset SO(\bot_p M )$ and we will show how to
construct the group in this paper. $\hat{G}_p \subset SO(\bot_p M )$   acts {\it polarly} on $\bot_p M$, this is to say
that there exists some linear subspace of $\bot_p M$ which meets every orbit and meets it perpendicularly
at every intersection point:

  {\bf Theorem}\hspace{0.3cm}
{\it 
Let $(G,K)$ be a Riemannian symmetric pair of compact type, such that 
$N=G/K$ is simply connected. Let $N$ be endowed with the    metric coming from the Killing form on $\mathfrak{g}$,
the Lie algebra of $G$. (This metric we will call ``Killing metric".)   Let  $\tilde{M} \subset N$ be an equifocal
manifold and
$M $ a focal manifold. Let $\hat {G}_p$ be as in Definition \ref{milhunzwan}, $\xi \in \bot_p M$ so that $\exp^N (\xi
)\in \tilde{M}$  and   $B_{\xi}:= \{ (\|^1_0 c) \hat{G}_p \xi |\
c:[0,1]\rightarrow M, c(0)=p  \}$. Then, we get $\tilde{B}_{\xi}:= \exp^N
(B_{\xi})=\tilde{M}$.}

  Moreover, we will find a characterization for the focal manifolds of equifocal submanifolds:

  {\bf Lemma}\hspace{0.3cm}  {\it 
Let $(G,K)$ be a Riemannian symmetric pair of compact type, such that 
$N=G/K$ is simply connected. Let $N$ be endowed with the Killing metric. Let  $\tilde{M} \subset N$ be an equifocal
manifold  and $M$ a focal manifold. Then,
$\hat{G}_p$ and the normal  parallel transport of $M$ preserve the focal structure of  $M$. }

  The converse is also true:

  {\bf Theorem}\hspace{0.3cm} {\it Let $M$ be simply connected and  $i:M\rightarrow N$ an isometric
immersion into a symmetric space  which has nonpositive or nonnegative sectional curvature. Let
$M$ admit an $\varepsilon$-tube for some $\varepsilon >0$. Let   $\xi$ (with $\| \xi \| <\varepsilon$) lie in
a principal orbit of $  \hat {G}_p$ and  
$\tilde{B}_{\xi}$ be the tube around $M$ defined as above.  Then, $\tilde{B}_{\xi}$ is an immersed
submanifold of $N$.
The immersion $\exp^N :B_{\xi} \rightarrow N$ has globally flat and abelian normal bundle. If $\hat{G}_p$ and
the normal parallel transport of
$M$ preserve the  focal structure of
$M$ then  the immersion is  equifocal and $M$ is a focal manifold of $\tilde{B}_{\xi}$.}

  In section 2 we will recall some facts about {\it  holonomy
systems} which were introduced by Simons in \cite{Sim}. In section 3 we
construct the Lie group $\hat{G}_p \subset \mbox{SO}(\bot_p M)$ and show
that it acts polarly on $\bot_p M$ by using the theory of holonomy
systems and algebraic curvature tensors. In section 4 we will prove   the second Theorem mentioned above. In section 5 we
consider a smooth  vector bundle
$E:[0,1] \rightarrow TN$ along some curve $c:[0,1]\rightarrow N$. Here, $N$ shall be of compact type.  If $\exp^N (E(t))$
is totally geodesic, flat and closed
in $N$ (and hence a flat torus) for all $t\in [0,1]$ then there exists
some curve $g:[0,1] \rightarrow I(M)$ into the isometry group of $N$  with $g(t)_{\ast}E(0)=E(t)$ for all $t$.
 This fact will be important in section 6, where we
will prove the first Theorem and the Lemma  mentioned above.

\section{Preliminaries}

   Let  $N$  be  an $n$-dimensional Riemannian  symmetric space, $M$ an $m$-dimensional  Riemannian
manifold   and   $i: M \rightarrow N$ an isometric immersion. We will identify $p\in M$ with $i(p) \in N$.
Throughout this paper
$\pi^{\bot} :\bot M \rightarrow M$ denotes  the  normal bundle of $M$ and  
$\nabla^{\bot}$ denotes the  covariant derivative of the normal bundle
induced  from the Levi Civita derivative of $N$.  
$(\|^1_0 c)$ denotes  the induced parallel transport along piecewise smooth curves  $c
:[ 0,1] \rightarrow M$ and $R^{\bot}$ the   normal curvature tensor. All curves are
piecewise smooth. Moreover,
$\phi_p := \{ (\|^1_0 \gamma ): \bot_p M \rightarrow \bot_p M
\ |\ \gamma (0)=\gamma (1) =p\}$ denotes  the   {\it normal holonomy group of $M$ in $p$} and $\phi^{\ast}_p$ it's
identity component, the so-called   {\it reduced  normal holonomy group of  $M$ in $p$}. Furthermore,  
$\pi : TN
\rightarrow N$ denotes the tangent bundle of 
$N$, $\nabla^N$ or $\nabla$ the  Levi Civita connection, $(|^1_0c
)$ the induced
parallel transport along curves  $c:[ 0,1] \rightarrow N$. As this parallel transport is different from the 
$\nabla^{\bot}$-transport (the so-called {\it normal  parallel transport}) we have chosen another symbol.
$R^N$ or $R$ denotes the curvature tensor of $\nabla^N$ and $A$ the shape operator of $M$ in $N$. For $\xi_p \in \bot_p
M$ we define  
\begin{eqnarray*} \widetilde{\mbox{\rm Hol}}_{\xi_p}(M) & := & \{(\|^1_0c)\xi_p\ |\ c:[0,
1]\rightarrow M,c(0)=p\} ,\\
 \mbox{\rm Hol}_{\xi_p}(M) & := & \exp^N(\widetilde{\mbox{\rm Hol}}_{\xi_ p}(M))
\end{eqnarray*}
 
   and will call both of them a {\it holonomy tube} around $M$.  

  In the following sections we will assume that   $M$ has an
{\it $\varepsilon$-tube}  in $N$  for some $\varepsilon >0$. This means that   for all $p \in M$  there exists an open
neighbourhood  $V$ of $p$ in $M$ such that
$\exp ^{N} | \{ v \in \bot M \ | \  \| v \| < \varepsilon , \pi^{\bot} (v) \in V
\}$ is a diffeomorphism onto an open subset $\tilde{V} \subset N$.

  Now, let $\eta \in \bot_p M$ be fixed. In the following we will identify $ T_{\eta}(\bot M) $ and
$ T_p M \oplus \bot_p M $ via the isomorphism $\dot{v}(0)    \mapsto  (\pi^{\bot}_{\ast} \dot{v}(0) , \nabla^{\bot}
_{\partial   }v(0) )$ for any curve $v:{\RR} \rightarrow \bot M$ such that $ v(0)=\eta$.

\begin{lemma}\label{drei} Let $\gamma$ be the geodesic in $N$  such that  $\dot{\gamma}(0)=\eta$. Let $R^N_{\eta}:T_p N
\rightarrow T_p N$ be the Jacobi operator in the direction of $\eta$. As $R_{\eta}^N$ is a selfadjoint operator we
can choose an orthonormal basis $\{ w_1 ,\ldots ,w_n \}$ of $T_pN$ consisting of eigenvectors of $R_{\eta}^N$. Let
$\lambda_1 ,\ldots ,\lambda_n$ be the corresponding eigenvalues. Then, 
\[   \exp^{\bot}_{\ast |
\eta} = ( |^1_0 \gamma )\circ (\bar{D}^{\eta} \oplus D^{\eta}):T_pM \oplus
\bot_pM \rightarrow T_{\exp^N (\eta )} N ,\]
  where \begin{eqnarray*} D^{\eta}: \bot_p M \rightarrow T_p N, & & \ z \mapsto \sum^n_{h=1} \sin
_{\lambda_h}(1)\langle z,w_h \rangle w_h \\
  \bar{D}^{\eta}: T_p M \rightarrow T_p N, & & \ z \mapsto \sum^n_{h=1} \cos
_{\lambda_h}(1)\langle z,w_h \rangle w_h - \sin_{\lambda_h}(1) \langle A_{\eta}(z),w_h  \rangle w_h.\\
\end{eqnarray*} 
  Here, $\cos_{\lambda}$ and $\sin_{\lambda}$ are the solutions
of the differential equation $y'' =-\lambda y$ with $\lambda \in {\RR}$, 
$\cos_{\lambda}(0)=1,\ \cos_{\lambda}' (0)=0 ,\ \sin_{\lambda}(0)=0 $ and
$\sin_{\lambda}' (0)=1$. 
\end{lemma}

  In the following sections we need to know some facts   about {\it holonomy systems} which were introduced by Simons
in \cite{Sim}. Let us therefore recall the definition and some of its properties: Let $V$ be a Euclidean vector space, $R$ an {\it algebraic curvature
tensor} on $V$ (this is to say a $(1,3)$-tensor on $V$ that satisfies the curvature equalities (including the first
Bianchi identity)). A compact subgroup  $G\subset SO(V)$ is called a {\em holonomy group}, if $
R_{x,y}\in\mathfrak{g}$ for all  $x,y\in V$. A triple   $S=[V,R,G
]$ consisting of a Euclidean vector space  $V$,
an algebraic curvature tensor $R$ and a connected holonomy group $G$ is called a
{\it  holonomy system}. The reduced holonomy group of a Riemannian manifold is a holonomy group in the sense of Simons
(this follows from  \cite{AS}).    $S$ is called  {\it symmetric} if  $g^{-1}R(g(),g())g=R$ for all
$g\in G$. For example, this is the case for the holonomy group of a locally symmetric space. 
$S$ is called {\it reducible} if $G$ acts  reducibly on $V$.
Every irreducible symmetric holonomy system corresponds to a simply connected, irreducible symmetric space
$N=\tilde{G}/H$ (where $\tilde{G} =I_0 (M)$, the identity component of the isometry group of $M$) in such a way that
  $T_p N \cong V$  (for  $p:= [H]\in N$), that  $R$ is the curvature tensor of  $N$ in $p$ and the action of  $G$ on  $V$
is the isotropy action of $H$ on
$T_{p}N$. For any  irreducible holonomy system $S=[V,R,G]$   whose scalar curvature $s(R)$ does not vanish 
there exists a curvature tensor $\tilde{R}$ such that $\tilde{S}=[V,\tilde{R},G]$
is symmetric:

\begin{lemma}[Simons]\label{hunvierzehn}
Let $S=[V,R,G]$ be an irreducible holonomy system, $\mathfrak{g}^R$ the Lie subalgebra of   $\frak{so}(V)$, which is
generated by   $\{ R_{x,y}|\ x,y\in V\} $ and $G^R$ the connected Lie subgroup of  $SO(V)$ whose Lie algebra is
$\mathfrak{g}^R$. We assume
$s(R)\not= 0$. Then, $G^R =G=N(G)$, where $N(G)$ is the normalizer of  $G$ in $SO(V)$. Moreover there exists another
curvature tensor $\tilde{R}$ on
$V$, such that $\tilde{S}=[V,\tilde{R},G]$ is an irreducible symmetric holonomy system. If 
$S$ is symmetric then $S=\tilde{S}$. (Here $\tilde{R}:= \int_{N(G)}  g^{-1}R (g(),g())g $ 
with Haar measure $\int_{N(G)}$.)
\end{lemma}

  Let $(G,K)$ be a Riemannian symmetric pair of compact type (or of noncompact type), $N:= G/K$ a corresponding
symmetric space and $\mathfrak{g}=\mathfrak{k}\oplus \mathfrak{p}$ the Cartan decomposition  corresponding to $K$. We
assume $K$ to be connected. Then,
$\mbox{Ad}_G (K)$ acts polarly on $\mathfrak{p}$ and this action is called {\it $s$-representation}. Dadok showed the converse (compare
\cite{Dad}): Let $H$ act irreducibly and polarly on the Euclidean space $V$. Then, there exists an $s$-representation  of
some Lie group
$K \hookrightarrow \mbox{SO}(V)$  on $V$ which has the same orbits as $H$. The principal orbits of $s$-representations are
{\it isoparametric}, this means they have globally flat normal bundle and constant focal distances. If we identify $T_p
N$ (with $p:= [K] \in N$)  and
$\mathfrak{p}$,  then we get: \vspace{-0.3cm}

\begin{eqnarray*} 
\left[ x,y\right]  \mbox{\hspace{1.6mm}}  & = & - R^N (x,y)
  \mbox{\    for all } x,y\in T_p N ,\\
\left[ A,B\right] & = & AB -BA   \mbox{\     for all } A,B \in \mathfrak{k},\\
\left[ A,x\right]  \mbox{\hspace{0.6mm}}  & =  & A(x)  \mbox{\         for all } A\in \mathfrak{k}, x\in T_p N,
\end{eqnarray*}
  where $[ , ]$ is the Lie bracket of $\mathfrak{g}$.
If  $N=G/K$ is   of compact type then the negative of the Killing form  $B$  
is a positive definite  inner product on
 $\mathfrak{g}$. Let $\langle ,\rangle$ be the inner product on $\mathfrak{p}=T_p N$ and without
loss of generality
$\langle ,\rangle =-B|\mathfrak{p}$.  Let $z_0 \in T_p N$ and $M=Kz_0 \subset T_p N$ be the orbit of
$K$ through $z_0$  and $K_{z_0}$ the corresponding isotropy group  in $z_0 \in T_p N$. 
 Let furthermore $\mathfrak{h}$ be the   Lie algebra of $K_{z_0}$ and $\mathfrak{m}$ the   orthogonal complement of
$\mathfrak{h}$ in $\mathfrak{k}$ corresponding to  $-B|\mathfrak{k}$.  Then, we get:
\[ T_{z_0}M =\{ [\eta ,z_0 ]|\ \eta \in \mathfrak{k} \} =\{ [\eta ,z_0 ]|\ \eta \in \mathfrak{m}\}.\]
  One can easily check  that
\[ \bot_{z_0}M =\{ \eta \in \mathfrak{p}|\  [\eta ,z_0 ]=0\} =\{ \eta \in T_p N | \ R^N (\eta ,z_0 )=0 \}
.\]
   The normal parallel transport of the orbit
$M$  along the curve $\exp (t m)z_0$, $m\in \mathfrak{m}$, is given by
$\exp (tm)_{\ast}$. 

  If $(G,K)$ is of noncompact type we can consider it's dual pair and get the same results.

\section{Construction of the polar  action}

  In this section $N$ shall be a symmetric space  whose  sectional curvature  is nonpositive or nonnegative.
Let $M$ be a simply connected manifold and $i:M\rightarrow N$ an isometric immersion admitting an $\varepsilon$-tube for
some $\varepsilon >0$. 

\begin{lemma}\label{som}
Let $N$ be a symmetric space whose sectional curvature is nonnegative or nonpositive. Let   
$\eta ,a
\in T_q N$ be chosen arbitrarily. Then, $\langle R^N (\eta ,a )\eta ,a \rangle =0$ implies $R^N (\eta  ,a )=0$.
\end{lemma}

  {\it proof}\hspace{0.3cm}   Let $T_q N =W_0 \oplus \ldots \oplus W_k$
be the  de Rham decomposition of $T_qN$ and  $R^N_q=\sum_{i=0}^kR^
i_q$ the corresponding decomposition of the
curvature tensor of $N$. Here, $W_0$ shall be the Euclidean component. For $\eta ,a \in T_q N$ there exist $\eta^i ,a^i \in
W_i$ such that $\eta =\sum_{i=0}^k \eta^i $ and
$a=\sum_{i=0}^ka^i$. By $0=\langle R^N_q(\eta ,a)\eta ,a\rangle =\sum_{i=0}^k\langle R_q^
i(\eta^i,a^i)\eta^i,a^i\rangle$ we get $\langle R_q^i (\eta^i ,a^i )\eta^i ,a^i \rangle =0$
 for all $i$, because we assumed   the sectional curvature to be nonnegative or nonpositive.
Hence, $R^i_q (\eta^i ,a^i)=0$ (the $R^i_q$'s for $i\geq  1$ correspond to irreducible symmetric spaces) and therefore
$R^N_q (\eta ,a)=0$.
\hfill
$q.e.d.$\\

   We define
\begin{equation} \tilde{R}_p :=\mbox{pr}_{\bot_p M} \circ (R^N_p |\bot_p M),\end{equation}
  where $\mbox{pr}_{\bot_p M}: T_p N \rightarrow \bot_p M$ is the orthogonal projection. It follows   $
\langle \tilde{R}_p (x,y)z,w \rangle = \langle  R^N_p
(x,y)z,w \rangle  $  for all  $x,y,z,w \in \bot_p M$ and
$\tilde{R}_p$ is  an algebraic curvature tensor in the sense of Simons.   Let $c:[0,1] \rightarrow M$
be a piecewise smooth curve with $c(0)=p$. Then,
\begin{equation} (\|^1_0 c)^{-1}  \tilde{R}_{c(1)}(\|^1_0 c)=:
\tau_c (\tilde{R}_{c(1)})\end{equation}
  is an  algebraic curvature tensor on  $\bot_p M$, too. Let $\hat{R}_p$
be the set of all such curvature tensors:
\begin{equation}
\hat{R}_p := \{ \tau_c (\tilde{R}_{c(1)}) |\ c:[0,1]\rightarrow M,\ c(0)=p \}.
\end{equation} 
  Let ${\cal L}_p$ be the  Lie 
subalgebra of $\frak{so}(\bot_p M)$ generated by  $\{ \tilde{R}(x,y)|\    \tilde{R}\in \hat{R}_p;\ x,y \in \bot_p M \}$
and   $G_p$ be the connected  Lie subgroup of  $SO(\bot_p M)$ whose Lie algebra is   ${\cal L}_p$.   Let $\bot_p
M=V_0^p
\oplus \ldots \oplus V^p_{k}$ be the decomposition of    $\bot_p M$ into $G_p$-invariant, for $i\geq 1$ irreducible
subspaces,   where
 $V^p_0$  is the subspace on which   $G_p$ acts trivially.  If we define
$\tilde{R}^i:=\mbox{pr}_i \circ (\tilde{R}|V^p_i )$ with  the orthogonal projection $\mbox{pr}_i:\bot_p
M\rightarrow V_i^p$ then $\tilde{R} (x ,y) =\sum_{i=0}^k  \tilde{R}^i (x_i  ,y_i  ) $, where 
$x_i ,y_i \in V_i^p$ such that $x =\sum_{i=0}^k x_i $ and
$y=\sum_{i=0}^k y_i$.  
 
  Now, let ${\cal L} ^{\tilde{R}^i}$ be the  Lie subalgebra of $\frak{so}(\bot_p M)$, which is generated by all
$\tilde{R}^i_{x_i,y_i}$ for $\tilde{R}\in \hat{R}_p$ and
$x_i,y_i\in V^p_i$.

\begin{lemma}\label{hunachzehn}
The spaces  ${\cal L}^{\tilde{R}^i},\ i=1,\ldots ,k$, are ideals in ${\cal L}_p$. 
If  $G_p^i $ is the connected   Lie 
subgroup   of $SO(\bot_p M)$ with  Lie algebra
${\cal L}^{\tilde{R}^i}$, then $ G_p=\{ \mbox{Id}_{V_0^p}\} \times G_p^1\times \ldots \times
G_p^{k}$ is a  direct  product, where
$G_p^i$ acts trivially  on $V^p_j$ for $i\not= j$ and  as an irreducible
$s$-representation on $V_i^p $ for
$i\geq 1$.
\end{lemma}

  {\it proof}\hspace{0.3cm} The proof of the fact that   $G_p=\{ \mbox{Id}_{V_0^p}\} \times G_p^1\times \ldots \times
G_p^{k}$ is the same as in  \cite{Sim} (compare the
Corollary on page 218 in \cite{Sim}).  And $G_p^i$ is compact for all $i=1,\ldots ,k$ by   \cite{KN1}, appendix (5).
For every $i\not= 0$ there exists a tensor $\tilde{R}\in \hat{R}_p$ such that
$\tilde{R}^i := \mbox{pr}_i \circ (\tilde{R}|V_i^p )\not= 0$. If $e_1 ,\ldots ,e_r$ is
an orthonormal  basis of $V_i^p$ then the
expression 
$\sum_{k,h=1}^r\langle\tilde {R}^i(e_h,e_k)e_h,e_k\rangle$ is a sum of nonnegative (or
nonpositive) numbers, where at least one number is   positive (negative) (compare 
Lemma
\ref{som}). So, the scalar curvature of  $\tilde {R}^i$ on $V_i^p$ is nontrivial. By Lemma \ref{hunvierzehn} there exists
some curvature tensor $\bar{R}^i$ such that  $ [V_i^p,\bar{R}^i,G_p^i]$ is a
 symmetric holonomy system. Thus,
 we see that $G_p= \{ \mbox{Id}_{V_0^p}\} \times  G_p^1 \times \ldots\times G^k_p$ acts as a product of
  $s$-representations and a trivial factor on  $V_0^p\oplus\ldots
\oplus V_k^p$.  \hfill $
q.e.d.$

\begin{lemma}\label{hunneunzehn}
$\phi_p^{\ast}$ leaves the subspaces $V_i^p$   invariant     and $(\phi^{\ast}_p |V^p_i )\subset G^i_p$ for
$i\geq 1$. Moreover,  there exists a decomposition $V_0^p =W_0 \oplus \ldots \oplus W_l$ into
$\phi_p^{\ast}$-invariant subspaces and a decomposition $\phi_p^{\ast}|V_0^p =\{ \mbox{Id}_{W_0} \} \times
H_1 \times \ldots \times H_l$ (direct product) such that $\phi_p^{\ast}|W_i =H_i$ acts  as
an irreducible $s$-representation for $i\geq 1$. 
\end{lemma}

  {\it proof}\hspace{0.3cm} Let $c:[0,1]\rightarrow M$ be
a piecewise smooth loop with starting point in    $p$ and
 $\tilde{R}\in \hat{R}_p$.
By the construction of $\hat{R}_p$  we get $\tau_c (\tilde{R} )\in \hat{R}_p
$ (for the definition of $\tau_c (\tilde{R} )$ compare (2)) and it follows
\[ (\|^1_0 c)^{-1}\tilde{R}(x,y)(\|^1_0 c)=\tau_c (\tilde{R}  )
((\|^1_0 c)^{-1}(x),(\|^1_0 c)^{-1}(y)) \in {\cal L}_p
,\]
 for every  $x,y \in \bot_p M$. Thus, $ (\|^1_0 c)^{-1}(G_p )(\|^1_0 c)\subset G_p$. So,
$\phi_p^{\ast}$ is contained in the normalizer of $G_p$ in $SO(\bot_p M)$.   As $\phi^{\ast}_p$ is connected  the 
  $V_i^p$'s are  $\phi_p^{\ast}$-invariant. So, $(\phi^{\ast}_p |V_i^p )$ is contained in the normalizer of
$G^i_p$ for $i\geq 1$ and by Lemma \ref{hunvierzehn} we get  $
(\phi^{\ast}_p|V^p_i )\subset G_p^i $ for
$i\geq 1$. 

   The construction of  $G_p$ described above can be done for any 
 point $p\in M$.  We define
\[ E:= \bigcup_{p\in M} V_0^p .\]
  By the construction of $G_p$ and of $V_0^p$ we get that $E$ is invariant under the $\nabla^{\bot}$-parallel
translation. Hence,  
$E$ is a  parallel smooth  subbundle of  $\bot M$.   In order to prove that  $\phi^{\ast}_p|V_0^p$ acts as a product of
$s$-representations and a   trivial factor we   define an algebraic curvature tensor  ${\cal
R}_p^{\bot }$ on $E_p$, such that
$S:=[E_p,{\cal R}_p^{\bot},\phi_p^{\ast}|V^p_0 ]$ is a  holonomy sy\-stem.    Let therefore $C^{\infty }(E)$ be the set of
all sections from $M$ to $E$. The curvature tensor
$R^{\bot}$ on $\bot M$  maps  $E_p$ to $E_p$  because $E$ is invariant under the $\nabla^{\bot}$-parallel
translation. Hence, we can define ${\cal R}^{\bot }:\ C^{\infty }(E)^3\rightarrow C^{\infty }(E)$ by
\[ {\cal R}_p^{\bot }(\xi_1,\xi_2)\xi_3= \sum^m_{j=1}R_p^{\bot
}(A_{\xi_1}(e_j),A_{\xi_2}(e_j))\xi_3\]
 {\rm for $p\in M,\ \xi_1,\xi_2,\xi_3 \in E_p$,
where $A$ is the shape operator of $M \subset N$  and $\{ e_1,\ldots ,e_m\}$ is any orthonormal basis of $T_pM$. 
As   $R^N(\xi ,\eta )=0$ for all  $\xi ,\eta\in E_p$ (by the construction of $E_p$), the Ricci equation has a very simple
form: 
\[ \langle R_p^{\bot}(x,y)\eta ,\xi \rangle =\langle [A_{\xi },A_{\eta}]x,y\rangle \]
  for $\eta ,\xi \in E_p ,x,y \in T_p M$.  Now, the proof that 
$\phi^{\ast}_p|V_0^p$ acts as a product of $s$-representations and a  
trivial factor is exactly the same as in \cite{Ol1}. \hfill $q.e.d.$

  {\bf Remark}\hspace{0.3cm} In general the Ricci equation is more complicated and therefore it was not possible
until now to generalize
 the results in \cite{HOT} and \cite{Ol1}  to arbitrary ambient spaces.

\begin{definition}\label{milhunzwan}{\rm Let
$\hat{G}_p \subset SO(\bot_p M)$ be the connected   Lie subgroup whose Lie algebra is generated by  
${\cal L}_p$ and the Lie algebra of  $\phi_p^{\ast}|V_0^p$. 
Then, $\hat {G}_p$ acts on  $V_1^p\oplus\ldots\oplus V^p_k$ as a product of
$s$-representations, because
$\hat{G}_p |V_1^p \oplus \ldots \oplus V_k^p =G_p$. And $\hat{G}_p$ acts on 
$V_0^p$ as a product of $s$-representations and a trivial factor. Because
$\hat{G}_p |V_0^p =\phi_p^{\ast}|V^p_0$.}
\end{definition}

  By construction, the    normal holonomy group  $\phi_p^{\ast}$ is  contained in   $\hat{G}_p$.
 For submanifolds of the Euclidean space   $\hat{G}_p$ and $\phi^{\ast}_p$ are the same. For submanifolds of the sphere
this is not necessarily true, because
$\hat{G}_p$ acts transitively on the unit sphere  in $\bot_p M$ but $\phi_p^{\ast}$ might not. This corresponds to the fact
that   isoparametric submanifolds of codimension $\geq 2$ of the sphere are not equifocal (because they do not have abelian
normal bundle).

\section{The tube $\tilde{B}_{\xi}$}

     $M$ shall be simply connected unless otherwise stated. At the end of this section we will think about the
case where $M$ is not simply connected (compare Corollary \ref{kyo}).

\subsection{Definition of the tube $\tilde{B}_{\xi}$}

   By the construction of $\hat{G}_p$ we get
\begin{equation} (\|^1_0 c) \hat{G}_p (\|^1_0 c)^{-1} = \hat{G}_{c(1)} \end{equation}
  for piecewise smooth curves
$c:[0,1]\rightarrow M ,c(0)=p$. Let  $\xi \in \bot_p M$ and $B_{\xi} \subset \bot
M$ be the subset   given by 
\begin{equation} B_{\xi} := \{ (\|^1_0 c) \hat{G}_p \xi |\ c:[0,1] \rightarrow M,c(0)=p
\}.\end{equation}
  If   $c:[0,1]\rightarrow M$ is any curve with  
 $c(0)=p,c(1)=q$ then by (4) we can conclude 
$(B_{\xi})_q=\hat {G}_q(\|^1_0c)\xi ,$   where $(B_{\xi})_q$ is the fibre over $
q$.  As $\hat {G}_p\xi$
is compact $B_{\xi} \subset \bot M$
is an embedded submanifold. If $M$ admits an $\varepsilon$-tube for some $\varepsilon 
>0$  and $\|\xi\|<\varepsilon$    then
$\exp^N |B_{\xi} :B_{\xi} \rightarrow N$ is an immersion and
\begin{equation} \tilde{B}_{\xi}:= \exp^N (B_{\xi})\end{equation}
  is an  immersed submanifold.  

  {\bf Remark}\hspace{0.3cm}  Let us endow $\bot M$ with the Sasaki metric (this means $\langle
X,Y\rangle_{\bot M} := \langle \pi^{\bot}_{\ast} (X),\pi^{\bot}_{\ast} (Y)\rangle_N + \langle K^{\bot} (X),K^{\bot}
(Y)\rangle_N$ for $X,Y \in T_{\eta}(\bot M), \ \eta \in \bot_q M$ and the connection map $K^{\bot} :T(\bot M)
\rightarrow \bot M$). By the fact that  $B_{\xi}$ contains the holonomy tube $\widetilde{\mbox{Hol}}_{\eta} (M)$ for any
$\eta
\in (B_{\xi})_q $ we get 
$\bot_{\eta} B_{\xi} \subset \bot_q M$, where  $\bot_{\eta} B_{\xi}$ is 
the normal space of  
$B_{\xi} \subset \bot M$.  By $\bot_{\eta}(\hat{G}_q \eta )$ we will denote  the normal space of  $\hat{G}_q \eta
\subset \bot_q M$ in $\bot_q M$.  Then, we get $\bot_{\eta} B_{\xi}= \bot_{\eta}(\hat{G}_q
\eta )$.

\subsection{Abelian normal bundle}

\begin{definition}
{\rm Let $\eta \in T_q N$. The} commutator of $\eta$ in $T_q N$  {\rm is the set of all $\nu
\in T_q N$ such that $R^N (\eta ,\nu )=0$. We say that a subspace   $A\subset T_q N$} is abelian {\rm if $R^N (\xi, \nu
)=0$ for all $\xi ,\nu \in A$}.
\end{definition}

  If $A\subset T_q N$ is abelian then $\exp^N (A)$ is a flat, totally geodesic submanifold in $N$, a so-called
{\it flat}.

\begin{lemma}\label{hundreizwanzig}
Let $\eta \in \bot_q M$ be chosen arbitrarily. 
Then, $R^N (\eta ,a )=0$ for all $a \in    \bot_{\eta} (\hat{G}_q \eta ) \subset
\bot_q M$ (compare the last Remark). Thus, the normal space of the orbit $\hat{G}_q \eta$  in $\bot_q M$ is contained in
the commutator of  $\eta $ in $T_q N$. The normal space is abelian if
$\eta$ lies in a principal orbit of
$\hat{G}_q$. 
\end{lemma}

  {\it proof}\hspace{0.3cm}  We fix
$\eta
\in
\bot_q M$ and $a  \in  \bot_{\eta} (\hat{G}_q \eta )$.    As $\tilde{R}_q (x,y)\eta$ is tangent to
$\hat{G}_q \eta$ by construction of $\hat{G}_q$ we get $ 0=\langle \tilde{R}_q (x,y)\eta ,a \rangle = \langle R^N_q
(x,y)\eta ,a \rangle $ for all  $x,y\in \bot_q M$, especially for $a,\eta$. As $N$
is a symmetric space with nonnegative or nonpositive sectional curvature we get
   $R^N_q(\eta ,a)=0$  (compare Lemma \ref{som}).  Thus,  
$\bot_{\eta} (\hat{G}_q \eta )$ is   contained in the commutator of $\eta $ in
$T_q N$.  This is true for any orbit, not only for principal ones. 

  Now, let $\eta \in \bot_q M$  be in a principal orbit of $\hat{G}_q$. We choose $\nu ,\mu
\in  \bot_{\eta}(\hat{G}_q \eta ) $. Especially, $\mu \in
 \bot_{\eta}(\hat{G}_q \eta )  \subset  \bot_{\nu}(\hat{G}_q \nu ) $ because
$\hat {G}_q$ acts polarly on $\bot_qM$ (compare Definition   \ref{milhunzwan}).
 So, we see $R^N (\nu ,\mu )=0$. \hfill
$q.e.d.$ 

  {\bf Remark}\hspace{0.3cm} Let $\eta\in\bot_qM\cap B_{\xi}$ (where $\xi$ does  not necessarily lie  in a
principal orbit)  and $A\subset  \bot_{\eta}(\hat{G}_q \eta )$ be an abelian subspace with
$\eta \in A$. Then:
\[ \exp_{\ast |\eta}^{\bot}  (0,A)= T_{\exp^N (\eta )}(\exp^N (A))\subset \bot_{\exp^N
(\eta )}(\tilde{B}_{\xi}),\]

  because  $(\bar{D}^{\eta} \oplus D^{\eta})(0,A) =0 $ and $(\bar{D}^{\eta} \oplus D^{\eta})(T_q
M,A^{\bot})\subset T_q M\oplus A^{\bot}$, where $A^{\bot}$ is the orthogonal complement of $A$ in $\bot_q M$ (compare 
Lemma
\ref{drei}).

  Now, let $\xi $   be in a principal orbit of   $\hat{G}_p$ and $\eta \in \bot_q M
\cap B_{\xi}$.  The normal space $A:= \bot_{\eta}(\hat {G}_q\eta 
)$ of the orbit $\hat {G}_q\eta$  in $\bot_qM$ therefore is abelian in 
$T_qN$ (compare Lemma \ref{hundreizwanzig}). By dimension arguments we get:
\begin{equation}  T_{\exp^N (\eta )}(\exp^N (A))= \bot_{\exp^N
(\eta )}(\tilde{B}_{\xi}) .\end{equation}

\begin{corollary}\label{key}
If  $\xi$ lies in a principal orbit of   $\hat {G}_p$ then the immersed  tube $
\tilde {B}_{\xi}$ has abelian normal bundle
 in $N$.
\end{corollary}

\subsection{Flat normal bundle}

  From now on $\xi\in\bot_pM$ shall lie in a principal orbit of $
\hat {G}_p$. Let $O(\bot M)$ be the bundle of orthonormal bases of $\bot M$. Let
$u=(p,u_1 ,\ldots ,u_l ) \in O(\bot M)$   and $B^u \subset O(\bot M)$   the subbundle  defined by
\[ B^u := \{ (\|^1_0 c) \hat{G}_p u|\ c(0)=p , c:[0,1] \rightarrow M \}.\]
   We identify $v=(q,v_1 ,\ldots ,v_l ) \in (B^u)_q$ with the isometry $\bot_p M \rightarrow \bot_q M ,\
(u_i \rightarrow v_i )_{i=1,\ldots ,l}$. By
$(\|^1_0 c) \hat{G}_p (\|^1_0 c)^{-1} = \hat{G}_{c(1)}$
  for all $c:[0,1]\rightarrow M ,c(0)=p$ (compare (4)),
 it follows that the fibre over $q\in M$ is equal to $
\hat {G}_q (\|^1_0 c) u$. The elements   in
$B^u$ leave  the tube  $B_{\xi}$ invariant and $B_{\xi}$ is induced by the ``action" of  $B^u$
on $\bot M$. To be more concrete: 
\begin{equation} B_{\xi}=\{ \psi \xi | \ \psi \in B^u\} .\end{equation}
   Therefore, we often will denote some element of
$B_{\xi}$ by $\psi \xi$.  In Lemma
\ref{hunsechszwanzig} we will see that the elements   in
$ B^u$ induce the  normal
parallel transport of $\tilde{B}_{\xi}$. And this helps us to show that 
$\tilde{B}_{\xi}$ has constant focal distances if and only if all elements in 
$B^u$ preserve the focal structure of   $M$   (compare Lemma
\ref{hundreidreissig}).

\begin{lemma}\label{hunsechszwanzig}
Let $\xi\in\bot_pM$   lie in a principal orbit of  $\hat {G}_p$ and let
 $\eta\in (B_{\xi})_q$. We choose $
\nu\in\bot_{\eta}(\hat{G}_q \eta ) \subset \bot_q M$.    Let $c:{\RR} \rightarrow M$ be a curve with
$c(0)=q$ and
$\psi :{\RR} \rightarrow B^u$ a lift with $\psi (0)=\mbox{Id}_{\bot_q M}$.  We claim that
\[ t\mapsto \exp^{\bot}_{\ast | \psi (t) \eta } (0, \psi (t) \nu ) \]
  is the   parallel  normal   field of $\tilde {B}_{
\xi}\subset N$ through
$\exp^{\bot}_{\ast | \eta }(0, \nu )$ along
$t\mapsto \exp^N (\psi (t)\eta )$. 
\end{lemma} 

  {\it proof} \hspace{0.3cm}  Let $w_1 ,\ldots ,w_k$    be a basis of the normal space 
$A:= \bot_{\eta}(\hat{G}_q \eta ) \subset \bot_q M$. We define  
\begin{eqnarray*} \tilde{z}(s,t) & := & \exp^{\bot}_{\ast | s \psi (t) \eta} (0, \psi (t) \nu ) ,\\
\tilde {w}_i(s,t) & := & \exp^{\bot}_{\ast | s \psi (t) \eta} (0, \psi (t) w_i  ),\ i=1,\ldots ,k .
\end{eqnarray*}

  Then, $s\mapsto\tilde {z}(s,t)$ and $s\mapsto \tilde{w}_i (s,t)$ are  parallel along the geodesic $
s\mapsto\exp^N(s\psi (t)\eta )$ for all $t$, because  $A$ is abelian and therefore $\exp^N (
A)$ a flat, totally geodesic submanifold.  The vectors
$\tilde{w_i}(1,t),i=1 ,\ldots k$, span the normal space
$\bot_{\exp^N (\psi (t)\eta )}(\tilde{B}_{\xi})$ (compare (7)).  In order to show that
$  \tilde{z}(1,t) = \exp^{\bot}_{\ast |\psi (t) \eta } (0, \psi (t) \nu )$
  is a  parallel normal vector field of $\tilde {B}_{\xi}$, it is sufficient to show    
\[\left\langle\left.\frac {\nabla^N}{dt}\right|_{t=w}\tilde z(1,t
),\tilde w_i(1,w)\right\rangle =0\]
  for $i=1,\ldots ,k$. For   $s=0$
we get:
\begin{eqnarray*}
& & \left\langle\left. \frac{\nabla^N}{dt}\right|_{t=w}\tilde{z}(0,t),\tilde{w_i}(0,w)\right\rangle \\
 & &  =  \left\langle\left. \frac{\nabla^N}{dt}\right|_{t=w}\exp^{\bot}_{\ast}
\left( \left. \frac{d}{dx}\right|_{x=0}(x\mapsto x\psi (t)     \nu )\right) ,\tilde{w}_i
(0,w)\right\rangle \\
 & & =  \left\langle\left. \frac{\nabla^N}{dx}\right|_{x=0}\exp^{\bot}_{\ast}
\left( \left. \frac{d}{dt}\right|_{t=w} (t\mapsto x\psi (t)\nu )\right) ,\tilde{w}_i
(0,w)\right\rangle .\\
\end{eqnarray*}

  But $t\mapsto \exp^{\bot} (x\psi (t)\nu )$ is a curve into  $\tilde{B}_{x\nu }$,
 the tube through $\exp^N (x\nu )$. In the   Remark in section 4.2 we saw that $ \exp^{\bot} (\psi
(w)A)$ meets the tube
$\tilde{B}_{x\nu }$ 
 perpendicularly. Hence,    $x\mapsto\exp^{\bot}_{\ast}(\frac d{d
t}|_{t=w}(x\psi (t)\nu )))$ is a normal vector field of   $\exp^{\bot} (\psi (w)A)$. As $\exp^{\bot} (\psi (w)A)$
is totally geodesic in $N$  (and hence its shape operator vanishes) we get that
$\left.  \frac{\nabla^N}{dx}\right|_{x=0}\exp^{\bot}_{\ast}
(\frac{d}{dt}|_{t=w}(x\psi (t)\nu ))$ is perpendicular to   $\psi (w) A$ and therefore
\[  \left\langle\left. \frac{\nabla^N}{dx}\right|_{x=0}\exp^{\bot}_{\ast}
\left( \left. \frac{d}{dt}\right|_{t=w}(t\mapsto x\psi (t)\nu )\right) ,\tilde{w}_i
(0,w)\right\rangle =0.\]

  It remains to show that $s\rightarrow \left\langle
\left. \frac{\nabla^N}{dt}\right|_{t=w}\tilde{z}(s,t),\tilde{w_i}(s,w)\right\rangle $ is constant. Let therefore
$\exp^N (s\psi (t)\eta )=: d(s,t)=:d^t (s)=:d_s (t)$.
\begin{align*}
& 
\left.
\frac{d}{ds}\right|_{s=x}\left\langle
\left. \frac{\nabla^N}{dt}\right|_{t=w}\tilde{z}(s,t),\tilde{w_i}(s,w)\right\rangle
 = 
\left\langle
\left.
\frac{\nabla^N}{ds}\right|_{s=x}\left.
\frac{\nabla^N}{dt}\right|_{t=w}\tilde{z}(s,t),\tilde{w_i}(x,w)\right\rangle ,\\
\intertext{because $s\mapsto \tilde{w_i}(s,w)$ is parallel along $d^w (s)$,}
 & = 
\left\langle \left.  \frac{\nabla^N}{dt}\right|_{t=w} \left. \frac{\nabla^N}{ds}\right|_{s=x}\tilde{z}(s,t),
\tilde{w_i}(x,w)\right\rangle 
 + \left\langle R^N (\dot{d}^w (x),\dot{d}_x(w))\tilde{z}(x,w),\tilde{w_i}(x,w)\right\rangle ,\\
& = \left\langle R^N (\tilde{z}(x,w),\tilde{w_i}(x,w))\dot{d}^w (x),\dot{d}_x (w) \right\rangle ,\\
\intertext{because $s\mapsto \tilde{z}(s,w)$ is parallel along $d^w (s)$,}
& =  0,
\end{align*}

  because $\tilde{z}(x,w)$ and $\tilde{w_i}(x,w)$ are contained in an abelian subspace.  
\hfill $q.e.d.$

\begin{corollary}\label{hundreissig}
$\tilde{B}_{\xi}$ has a  globally flat normal bundle. To be more concrete: Let $\nu \in \bot_{\xi}(\hat{G}_p \xi
) \subset \bot_p M$. Then, $B_{\xi} \rightarrow \bot \tilde{B}_{\xi},\  \psi \xi \rightarrow \exp^{\bot}_{\ast | \psi  \xi
} (0, \psi  \nu ) $ is the parallel normal vector field of $\tilde{B}_{\xi}$ through $\exp^{\bot}_{\ast |    \xi }
(0,   \nu ) $.
\end{corollary}

  {\it proof}\hspace{0.3cm}  We have to show that $\psi \xi  \rightarrow \exp^{\bot}_{\ast | \psi \xi}(
0,\psi \nu )$  is   well defined.
  Let therefore $\psi$ and $ \varphi \in  B^u$ such that $\psi \xi =\varphi \xi$. Then, $\varphi^{-1}\circ \psi   \in
\hat{G}_p$. As $\xi$ lies in a principal orbit of  $
\hat {G}_p$
 the isotropy group $(\hat{G}_p)_{\xi}$ of $\xi$ acts trivially on the normal space of the orbit
 $\hat{G}_p \xi$ in $\bot_p M$ by what we can conclude  
$\varphi^{-1} \circ \psi   (\nu )= \nu$.  Moreover, $\psi \xi  \rightarrow \exp^{\bot}_{\ast | \psi \xi}(
0,\psi \nu )$ is parallel by Lemma \ref{hunsechszwanzig}. \hfill  $q.e.d.$

\subsection{Constant focal distances}

\begin{lemma}\label{huneindreissig}
Let $\rho\in\bot_qM$ be chosen arbitrarily (not necessarily in a principal orbit of $\hat{G}_q$). Then, {\rm
\[
\mbox{ker} (\exp^{\bot} _{\ast |
\rho})=\mbox{ker}(\exp^N_{\ast} |T_{\rho } B_{\rho} ) ,\] }
  where $T_{\rho}B_{\rho}$ is the tangent space of   $B_{\rho}\subset \bot M$.
\end{lemma}

  {\it proof}\hspace{0.3cm}  The   normal space
$\bot_{\rho}(B_{\rho})\subset\bot M$ is contained in the commutator of $\rho$ in $T_qN$   (compare Lemma
\ref{hundreizwanzig}). Thus, $R_{\rho}| \bot_{\rho}(B_{\rho}) =0$. Now,  the proof   follows   by
Lemma \ref{drei}.\hfill 
$q.e.d.$

\begin{definition}\label{un}
{\rm We say that a   linear isometry $ \psi :\bot_p M \rightarrow \bot_q M$} preserves the focal structure of $M$
{\rm  if for every focal vector  $\eta \in \bot_p M$ of multiplicity $m$  the vector $ \psi (\eta )$
is a focal vector of multiplicity  $m$, too.}
\end{definition}

\begin{lemma}\label{hundreidreissig}
The tube  $\tilde{B}_{\xi}$ has constant focal distances if and only if all elements   in $B^u$ preserve the
focal structure of $M$. This is equivalent to the fact that $\hat{G}_p$ and the normal parallel transport of $M$ preserve
the focal structure of  $M$.
\end{lemma} 

  {\it proof}\hspace{0.3cm}  Let $q\in M$ be chosen arbitrarily and   $\eta \in (B_{\xi})_q$. Let   $\nu \in
\bot_{\eta}(\hat{G}_q \eta ) \subset \bot_q M$ and $\tilde{\nu}:= \exp^{\bot}_{\ast | \eta }(0,\nu )\in \bot_{\exp^N
(\eta)}(\tilde{B}_{\xi})$. Let
$\rho :=\nu +\eta\in\bot_qM$.   As $\rho$ 
does not necessarily  lie in a principal orbit of $\hat {G}_q$,  the tube $
B_{\rho}$ possibly  has a smaller
dimension than  $B_{\xi}$. We can define a projection   $\Omega_{\xi ,\rho} :B_{\xi} \rightarrow B_{\rho}$ as follows:

  Let us choose    $\tilde{u}=(q,\tilde{u}_1 ,\ldots ,\tilde{u}_l ) \in B^u$ and   identify the elements   $(z,v_1
,\ldots ,v_l )\in B^u$ with the isometries $\bot_q M \rightarrow \bot_z M ,\ (\tilde{u}_i \rightarrow v_i )_{i=1,\ldots
,l}$.   Let \vspace{-0.5cm}

\[ \Omega_{\xi ,\rho}(\psi \eta ):=\psi \rho\]
  for all $\psi \in B^u$. This is well defined as one can see as follows:
  Let us choose $\psi_1$ and $\psi_2 \in B^u $ with $\psi_1
\eta =\psi_2 \eta$. It follows   $\psi_1^{-1}\psi_2  \in (\hat{G}_q )_{\eta}$. As 
$\hat {G}_q$ acts polarly on $\bot_qM$ and $\eta$ lies in a principal orbit of this polar action and  $
\rho\in \bot_{\eta}(\hat {G}_q\eta )  $ we get $
(\hat {G}_q)_{\eta}\subset (\hat {G}_q)_{\rho}$ (compare
\cite{PT}, page 81).   And so we get $\psi_1\rho =\psi_2\rho$, by what we can conclude that $
\Omega_{\xi ,\rho}$ 
is well defined. One can easily check that
$\Omega_{\xi ,\rho}$ is a   surjective submersion.

  Now, we want to assume that the elements in $B^u$ preserve the focal 
structure of $M$. Let $\hat{\nu}$ be the parallel  normal vector field of $\tilde { B}_{\xi}$
through $\tilde{\nu}$. In Corollary \ref{hundreissig} we saw   $ \hat{\nu}( \psi \eta )=
\exp^{\bot}_{\ast | \psi  \eta } (0, \psi  \nu ) $.  

  As $ \exp^N ( \bot_{\eta}(\hat{G}_q \eta ))$ is    flat and totally geodesic in $N$   
 and $\nu \in \bot_{\eta}(\hat{G}_q \eta
)$ we get   
\[ \exp^N (\rho ) =\exp^N (\eta + \nu  )= \exp^N (\exp^{\bot}_{\ast | \eta }(0,\nu ))= \exp^N(\tilde{\nu }) =\exp^N
(\hat{\nu}(\eta )).\]
  In the same way we get for all $\psi \eta \in B_{\xi}$:
\begin{equation} \exp^N (\Omega_{\xi ,\rho}(\psi \eta ))= \exp^N  (\psi \rho )
= \exp^N(\hat{\nu }  (\psi \eta )).\end{equation}

   So, we get  
\begin{align*}
 &  \mbox{dim}( \mbox{ker}(\exp^N_{\ast} |T_{\hat{\nu}(\psi \eta )}(\bot \tilde{B}_{\xi})))  =   \mbox{dim}(
\mbox{ker}(\exp^{\bot}_{\ast | \hat{\nu}(\psi \eta )} |(T_{\exp^N (\psi \eta )} \tilde{B}_{\xi}, 0))) \\
\intertext{because $\tilde {B}_{\xi}$ has abelian normal bundle. By (9) we get}
 & =   \mbox{dim(ker}(\Omega_{\xi ,\rho})_{\ast |\psi \eta }) + \mbox{dim(ker}(\exp^{N}_{\ast}|T_{\psi \rho} B_{\rho}))\\
 & =   \mbox{dim(ker}(\Omega_{\xi ,\rho})_{\ast |\psi \eta}) + \mbox{dim(ker}(\exp^{\bot}_{\ast |\psi \rho} ))
\end{align*}
  by Lemma \ref{huneindreissig}. As $\Omega_{\xi, \eta}$ is a surjective submersion and therefore has constant
rank, we get that  $\tilde{B}_{\xi}$ has constant focal distances if the elements in 
 $B^u$ preserve the focal structure of $M$.

  Now, we assume that $\tilde{B}_{\xi}$ has constant focal distances and
we want to prove that the elements  in $B^u$ preserve the focal structure of  $M$. 
Let us fix $\rho \in \bot_q M$   and  $\eta \in (B_{\xi})_q$. As  $\bot_q M = \bigcup_{\tilde{\eta }\in \hat{G}_q
\eta} \bot_{\tilde{\eta}}(\hat{G}_q
\eta ) $ we can find some  $\tilde{\eta} \in \hat{G}_q \eta$ such that $\rho \in
 \bot_{\tilde{\eta}}(\hat{G}_q \eta ) $. Let $\nu := \rho -\tilde{\eta}$ and $\tilde{\nu}:=
\exp^{\bot}_{\ast | \tilde{\eta}}  (0,\nu )$. The rest 
of the proof is the same as above.
\hfill
$q.e.d.$

\begin{theorem}\label{milvierhun}
Let $M$ be simply connected and  $i:M\rightarrow N$ an isometric immersion into a symmetric space  which has nonpositive
or nonnegative sectional curvature. Let
$M$ admit an $\varepsilon$-tube for some $\varepsilon >0$. Let furthermore
$\hat {G}_p\subset SO(\bot_pM)$ be the Lie group introduced in Definition \ref{milhunzwan},
$\xi\in\bot_pM$  with $\|\xi\|<\varepsilon$ so that $\xi$ lies in a principal orbit of $
\hat {G}_p$ and  
$\tilde{B}_{\xi}$ the tube around $M$ as in  (6).  Then, $\tilde{B}_{\xi}$ is an immersed
submanifold of $N$. The immersion $\exp^N :B_{\xi} \rightarrow N$ 
 has globally flat and abelian normal bundle. If $\hat{G}_p$ and the normal parallel transport of $M$ preserve the focal
structure of
$M$ then  the immersion is equifocal.
\end{theorem}

  {\bf Remark}\hspace{0.3cm} Let $M$ be as in Theorem
\ref{milvierhun} and $\psi \xi \rightarrow \exp^{\bot}_{\ast |\psi \xi} (0, -\psi \xi
)$ the parallel normal vector field of $\tilde{B}_{\xi}$ through $\exp^{\bot}_{\ast | \xi} (0, - \xi )$. Then,     
$\exp^N(\tilde{\xi }(\tilde {B}_{\xi}))= M$. Hence,
  $M$ is a focal manifold of the immersed equifocal manifold $\tilde{B}_{\xi}$ in $N$.

  {\bf Remark}\hspace{0.3cm} The assumption that  $\hat{G}_p$ and the  normal  parallel transport of $M$ both shall
preserve the focal structure of $M$ is necessary. This means it cannot be replaced by the weaker condition of ``constant
focal distances" in the sense of section 1.  Let 
$M\subset S^n$ for example be an  isoparametric submanifold with codimension $\geq 2$.   The   holonomy tubes
of
$M$ have at least the same codimension.   Therefore they cannot have abelian normal bundle and so they cannot be
equifocal. But the tube defined in this section has codimension 1. In general the group
$\hat{G}_p$ will not preserve the focal structure of $M$. Thus, the tube $\tilde{B}_{\xi}$ will not be equifocal in $N$.

  {\bf Remark}\hspace{0.3cm} The sectional curvature of  $N$ is bounded and one can apply the
Comparison Theorem of Rauch (compare
\cite{CE}, page 29) in order to show that  $M$ is complete if and only if the tube  $\tilde{B}_{\xi} \subset N$ with
the by $N$ induced metric is complete. The proof is similar to the one in \cite{HOT}.

  {\bf Remark}\hspace{0.3cm} Until now we assumed    $M$ to be simply connected. If this is not the case
then let $f:\tilde {M}\rightarrow M$ be the universal covering of $
M$. It is clear that  $i\circ f:\tilde {M}\rightarrow N$ is an 
  isometric immersion  which admits an   $\varepsilon$-tube in $N$ if and only if 
  $i$ admits one. Now, we fix $\tilde {p}\in\tilde {M}$ and $p:=f
(\tilde {p})\in M$. We   identify  
$\tilde{p}$ and $p$ with $i(f(\tilde{p}))=i(p)\in N$. Let   $\xi \in \bot_p M \cong \bot_{\tilde{p}}\tilde{M} \subset
T_{i(p)}N$.  Let  $B_{\xi} \subset \bot \tilde{M}$ be the tube around   $\tilde{M}$ in the sense of section 4.1. We  
define    $D_{\xi} \subset \bot M$ as follows:
\[ D_{\xi}:= \{ (\|^1_0 c) \hat{G}_p \xi |\  c:[0,1] \rightarrow M ,c(0)=p\} ,\]
  where  $\hat {G}_p$ is the connected Lie group introduced in Definition \ref{milhunzwan}.
The fibre of $D_{\xi}\rightarrow M$
over $p$ is not necessarily connected. But by construction it is clear that
\[ \exp^N (B_{\xi}) =\exp^N (D_{\xi}) .\]
  $\hat{G}_p$ and the normal parallel transport of  $M$    preserve the focal structure of
  $M$ if and only if the same is true for  $\tilde{M}$.  

\begin{corollary}\label{kyo}
Let      $\hat{G}_p$ be the connected Lie subgroup of
$SO(\bot_p M)$  generated by $G_p$ and the reduced normal holonomy group $\phi_p^{\ast}$ of $M$. Moreover, let $\xi\in\bot_
pM$ lie in a principal orbit of
$\hat{G}_p$ with $\| \xi \| <\varepsilon$ and 
$\tilde{B}_{\xi}:=\exp^N (D_{\xi})$  where $D_{\xi}\subset \bot M$ is defined as above. $\tilde{B}_{\xi}$ is  an immersed
submanifold of $N$.  The immersion $\exp^N :D_{\xi} \rightarrow N$  has   flat and abelian normal bundle. If
$\hat{G}_p$ and the normal  parallel transport of
$M$ preserve the focal structure of   $M$ then the immersion  is locally equifocal and $M$ is a focal manifold of
$\tilde{B}_{\xi}$.
\end{corollary}

\section{A smooth variation  of tori}

  Let $(G,K)$ be a  Riemannian symmetric pair of   compact  type and  $N=G/K$ a corresponding symmetric space.
Let $c:[0,1] \rightarrow N$ be a smooth curve and  $E(t) \subset T_{c(t)}N$ for all $t\in [0,1]$ an abelian subspace of
dimension  $d$ such that 
\[E:=\bigcup_{t\in [0,1]}E(t)\]
 is a  smooth vector bundle  along $c$. This   means there exists some smooth orthonormal basis vector field along 
$c$ that spans   $E(t)$ for all $t\in [0,1]$. Let us  fix   such a    smooth orthonormal basis 
vector field $\{v_1,\ldots ,v_d\}$.
 Let $\psi_t^{t_0}:E(t_0)\rightarrow E(t)$ be the isometry given by
$(v_i (t_0 ) \mapsto v_i (t))_{i=1,\ldots ,d}$.

\begin{theorem}\label{vierzehn}
If $\exp^N (E(t))$ for all $t\in [0,1]$ is a compact torus then there exists a lift $\tilde{b}:[0,1]\rightarrow G$ of
$c:[0,1]\rightarrow N$ such that $ \tilde {b}(t)_{\ast}E(0)=E(t).$  Hence, $\tilde{b}(t)$  maps $\exp^N (E(0))$ to
$\exp^N(E(t))$.
\end{theorem}

  {\bf Remark}\hspace{0.3cm} In   Theorem \ref{vierzehn} we can choose $
\tilde {b}$ such that
  $\tilde {b}(0)=e$ because if $\tilde {b}$ satisfies the assumptions of the Theorem
 then $\tilde{b}\cdot (\tilde{b}(0)^{-1})$ does so, too.  

   Let $G_{t}\subset G$   be the identity component  of the isotropy group  in
$c(t)$. Then, $G_{t}$ acts linearly on $T_{c(t)}N$ by $(z ,k)\mapsto k_{\ast}z $. We will usually write $kz$ instead of
$k_{\ast}z$. Let $z_0 \in T_{c(t)}N$ and  $(G_{t})_{z_0}:=\{ k\in G_{t} |\ kz_0 =z_0 \}$.

\begin{lemma}\label{siebzehn}
 The set $ \hat{W}_l :=  \{ z_0 \in E(t) | t\in [0,1], \dim ((G_{t})_{z_0} )\leq l\} $  is an open subset of $E$ for all
$l$. 
\end{lemma}

  {\it proof}\hspace{0.3cm}  We define $W_l := \{ z_0 \in T_{c(0)}N| \ \mbox{dim}((G_{0})_{z_0}) \leq l\}$. By
\cite{PT}, page 81, $W_l$ is open in   $T_{c(0)}N$.   Moreover,  $W_0 \subset W_1
\subset \ldots $. Let $\tilde{c}:[0,1] \rightarrow G$ be a lift of $c$ with $\tilde{c}(0)=e$.
Now, the proof  follows easily by the fact that $\{ z_0 \in T_{c(t)}N |\ \dim ((G_{t})_{z_0}) \leq
l \} = \tilde{c}(t)_{\ast} (W_l  )$.
  \hfill $q.e.d.$\\

  Let   $k:=\min\{l\in{\NN}_0|\ \hat {W}_l\not =\emptyset \}.$ Let furthermore $
t_0\in [0,1]$ and $\xi\in E(t_0)\cap\hat {W}_k$   be fixed, let $K:= G_{t_0}$ be the isotropy group in $c(t_0 )$ and  
\[\xi (t):=\psi^{t_0}_t(\xi ).\]
    Let $I_0 \subset [0,1]$ be  an open interval such that $\xi (t)\in \hat{W}_k$ for
all $t\in I_0$.

  Let $\hat{c}: I_0 \rightarrow G$ be a  smooth lift of $c|I_0$ with $\hat{c}(t_0 )=e$ such that
  $\hat{c}(t)c(t_0 )=c(t)$.
Then, $(t\mapsto \hat{c}(t)_{\ast}^{-1}(E(t)))=: \hat{E}(t)$  is a smooth vector bundle consisting of abelian subspaces of
$T_{c(t_0)}N$. We define
\[ h(t):= \hat{c}(t)
_{\ast}^{-1}(\xi
(t)).\]
   Let   $K \xi \subset T_{c(t_0 )}N$ be the orbit of $K$ through $\xi$. There exists an $r$-tube around
$K\xi$ in $T_{c(t_0 )}N$ for some $r\in {\RR}^{+}$. For $x\in K\xi$ let $S_x := \{\mu\in\bot_x(K\xi )|\ \|\mu
\|<r\}$.   As  $h$ is smooth and   $h(t_0 )=\xi$ we can find an open interval $J_0
\subset I_0$ such that $h(  J_0 )$ is contained in the $r$-tube around $K\xi$. Thus, we can find a smooth curve
$\gamma :J_0\rightarrow K\xi$ such that $h(t)\in S_{\gamma (t)}$ for all 
$t\in J_0$. Let  $o(t), t\in J_0$, be a smooth curve in $K$ with
$ o(t)\xi  =\gamma (t)$.

\begin{lemma}\label{zwanzig} 
$ \{
o(t)^{-1}(\hat{E}(t))|\ t\in J_0 \}$ is contained in an abelian subspace
of $T_{c(t_0 )}N$.
\end{lemma}

  {\it proof}\hspace{0.3cm}  $K_y \subset K_x$ for all $y\in
S_x$ and for all $x\in K\xi$ (compare \cite{PT}, page 81).  For $y\in S_{\gamma (t)}\cap\hat {E}(t)$ we get 
$\dim(K_y)\geq k=\dim(K_{\xi})=\dim(K_{\gamma (t)})$ (compare the definition of $k$). Therefore $K_y=K_{\gamma
(t)}$. Let
\begin{equation*}
 A(t) := \bot_{\gamma (t)}(K\xi )   = \{ \eta \in T_{c(t_0 )}N| R^N (\eta , \gamma (t))=0\}
\end{equation*}

  be the normal space of the orbit $K\xi$ 
 in $\gamma (t)$. Then,
\[ A(t)=\bot_y (Ky)=\{ \eta \in T_{c(t_0 )}N| R^N (\eta , y)=0\} \]
  for all $ y\in S_{\gamma (t)}\cap\hat {E}(t)$  (the proof is similar to the one that $K$ acts polarly on
$T_{c(t_0)}N$).  Especially, we get $\hat{E}(t)\subset
\{ \eta \in T_{c(t_0 )}N|\ R^N (\eta ,h(t))=0 \} = A(t)$ for all $t\in J_0$, because $h(t) \in  S_{\gamma
(t)}\cap \hat{E}(t)$ and
$\hat{E}(t)$ is abelian.
$S_{\gamma (t)}$ is open in 
$A(t)$ and so $S_{\gamma (t)}\cap
\hat{E}(t)$  is open in $\hat{E}(t)$. 
We can choose a basis $\{ v_1 , \ldots ,v_d \}$ of $\hat{E}(t)$ such that $v_i \in S_{\gamma
(t)}$ for all $i=1, \ldots ,d$. Thus,  we get 
\[  R^N (\hat{ E}(t),A(t))=0 .\]
  It  follows that $R^N (o(t)^{-1}(\hat{E}(t)),o(t)^{-1}(A(t)))=R^N (o(t)^{-1}(\hat{E}(t)),A(0))=0$.
Especially $ R^N (o(t)^{-1}(\hat{E}(t)), o(\tilde{t})^{-1}(\hat{E}(\tilde{t})))=0 $ for all $t,\tilde{t}
\in J_0$. Thus,  
$\{ o(t)^{-1}(\hat{E}(t))|\ t\in J_0 \}$ is contained in a maximal abelian subspace   of $T_{c(t_0 )}N$.\hfill
$q.e.d.$

  We denote by $B\subset T_{c(t_0 )}N$ the maximal abelian subspace. As $N$ is a symmetric space of compact type
$\exp^N(B)$ is compact and therefore a torus. By assumption  
\[\exp^N(o(t)^{-1}(\hat {E}(t)))=o(t)^{-1}\hat {c}(t)^{-1}(\exp^N
(E(t)))\]
   are compact tori, too.  So,
it only remains to prove:

\begin{lemma}\label{einzwanzig}
Let $T$ be a torus of dimension $l$, $T_0 T$ the tangent space  in $0\in
T$ and
$t\mapsto E(t)\subset T_0 T$ ($t\in J_0$) a smooth vector bundle of subspaces of
dimension $d\leq l$ such that $\exp^T (E(t))$ are subtori for all $t\in J_0$ . Then, $E(t)=E(0)$ for all $t\in
J_0$.
\end{lemma} 

  {\it proof}\hspace{0.3cm}   Let $\{ x_1 ,\ldots ,x_l \} $ be the lattice of $T$. Let
$G^{{\RR}}(d,l)$ be the
Gra{\ss}mannian manifold of all  $d$-planes in ${\RR}^l$. The subset $M \subset G^{{\RR}}(d,l)$   consisting of all
$d$-planes  spanned by basis vectors of the form   
\[a^1:=\sum_{i=1}^la_i^1x_i,\ldots ,a^d:=\sum_{i=1}^la_i^dx_i\mbox{ \rm with }
a_i^1,\ldots ,a_i^d\in\ZZ,\]
 
  (this means which contain a full lattice)  is countable.
$E:J_0
\rightarrow G^{{\RR}}(d,l)$ is smooth and $E(J_0 )\subset M$. Hence, we get  $E(t)=E(0)$ for all $t\in
J_0
$. \hfill $q.e.d.$

\begin{corollary}\label{achzehn}
Let $t_0 \in [0,1]$ such that $\hat{W}_k \cap E(t_0 ) \not= \emptyset$. There exists an open interval $ J_0 \subset [0,1]$
with $t_0 \in J_0$ and a smooth curve $w:J_0 \rightarrow G$ such that  $w(t)_{\ast}E(t_0 )=E(t  )$.
\end{corollary}

  {\it proof}\hspace{0.3cm}  $\hat{c}(t)o(t)=:w(t)$ is the curve we are looking for. \hfill $q.e.d.$

  {\it proof of Theorem  \ref{vierzehn}}\hspace{0.3cm}
  Let  $G_d^{\RR}(N)$ be the  Gra{\ss}mannian manifold of all $d$-plains in
the tangent bundle of $N$. Then, $G$ acts on   $G_d^{\RR}(N)$ by
\[ (g,\mbox{span}
\{ v_1 ,\ldots ,v_d \} )\mapsto \mbox{span}(\{ g_{\ast}v_1 ,\ldots ,g_{\ast} v_ d \} ) ,\]
   where $v_1 ,\ldots ,v_d$ are linearly  independent. As $G$
is compact the orbits are embedded submanifolds of  $G_d^{\RR}(N)$. 
And the projection of   $G$ onto any  orbit is a Riemannian submersion. Let $J_0 \subset [0,1]$ be the maximal interval
with $t_0 \in J_0$ such that $E(J_0 ) \subset G\cdot E(t_0 )$, where $G\cdot E(t_0 )$ is the orbit of $G$ through $E(t_0
)$.  By Corollary \ref{achzehn} we get that $J_0 $ is open (because the Corollary is true for any $t  \in [0,1]$ such that
$\hat{W}_k \cap E(t  )\not= \emptyset$). But $ t\mapsto E(t) $
 is a smooth curve into $G_d^{\RR} (N)$  by assumption ($E$ is a smooth vector bundle). Hence, $J_0$ is closed, and
therefore $J_0 =[0,1]$.  
As $t\mapsto E(t)$ is contained  in an orbit it can be liftet to a smooth curve   $\tilde{b}$ into $G$. 
\hfill
$q.e.d.$

\begin{lemma}\label{sechzig} 
Let $g: {\RR}  \rightarrow G$ be a  smooth curve with $g(0)=e$ and $E  \subset T_p N$ be an abelian subspace  such
that $\exp^N (E)$ is compact. Let furthermore $\frac{d}{dt}|_{t=0}g(t)p$ be perpendicular to $E$. 
 Then, the vector field defined by 
\[X:\exp^N(E)\rightarrow TN,\ y\mapsto\left.\frac d{dt}\right|_{t
=0}g(t)y\]
  is a normal vector field of $\exp^N(E)$.
\end{lemma}

  {\it proof}\hspace{0.3cm}  Let  $\gamma :\RR
\rightarrow\exp^N(E)$ be any geodesic. Then, $t\mapsto g(t)(\gamma 
)$ 
is a smooth variation of geodesics in $N$. Hence,  $X\circ\gamma$ is a Jacobi vector field. The orthogonal projection
$(X)^{\bot}\circ\gamma$ of $X\circ\gamma$ to $\exp^N(E)$ is a  Jacobi vector field, too. We especially get
$\langle X^{\bot}(\gamma (t)),\dot{\gamma}(t)\rangle =at+b$  with
$a,b
\in {\RR}$. As $\exp^N (E)$ is compact
  the norm of the vector field   $X$ is bounded. Thus, we see
$\langle X^{\bot}(\gamma (t)),\dot{\gamma}(t)\rangle =b$, by which $\langle (\nabla^N X 
^{\bot})(\dot{\gamma }(t)),\dot{\gamma}(t)\rangle =0$ follows. As this is true for any geodesic,
  $X^{\bot}$ is a   Killing vector field along $\exp^N(E)$
and therefore  trivial, because $\exp^N (E)$ is a flat torus and $0=(\frac{d}{dt}|_{t=0}g(t)p)^{\bot}=
X^{\bot}(p)$
by assumption. Hence, $X$ is a normal vector field along $\exp^N(
E)$.  \hfill $q.e.d.$

\begin{corollary} \label{milzehn}
Let $M$ be a submanifold of a  symmetric space $N$ of compact type and let $M$
have abelian and compact normal bundle. Here, we say that $M$ has   compact normal bundle  if $\exp^N (\bot_p M) $ is
compact for all $p\in M$.  Let
$c:{\RR} \rightarrow M$ be a curve      and
$g:{\RR} \rightarrow G$, $g(0)=e$ (where $G=I_0 (M)$), a curve which induces isometries of the tori $\exp^N(\bot_{c(0)}M
)\rightarrow\exp^N(\bot_{c(t)}M)$ (as in 
Theorem \ref{vierzehn}). These isometries induce the $\nabla^{\bot}$-parallel transport.  We
especially see that
$M$ has flat normal bundle because  $\phi^{\ast}_p$ ( $p\in M$) is contained in the isotropy group
 of the isometry group of a torus and therefore trivial. If $M$ is not simply connected then 
$M$ does not necessarily have globally flat normal bundle.
\end{corollary} 

  {\it proof}\hspace{0.3cm} Let $s\in {\RR}$, $
v\in\bot_{c(0)}M$ and    $\gamma :\RR\rightarrow\exp^N(\bot_{c(s)}
M)$ the geodesic with $\dot{\gamma }(0)=g(s)_{\ast}v$. Let $X$ be
 the vector field induced by $g$ along 
$\exp^N(\bot_{c(s)}M)$. That is to say  $X(q):=\frac d{dt}|_{t=s}
g(t)g(s)^{-1}q$ for all $q\in\exp^N(\bot_{c(s)}M)$.  Then,

\begin{eqnarray*}
\left. \frac{\nabla^N}{dt}\right|_{t=s}g(t)_{\ast}  v  
& = & \left. \frac{\nabla^N}{dt}\right|_{t=s}g(t)_{\ast}g(s)^{-1}_{\ast}g(s)_{\ast}  v  \\
& = & \left. \frac{\nabla^N}{dt}\right|_{t=s}
\left. \frac{d}{dw}\right|_{w=0}g(t)g(s)^{-1}\gamma (w)\\ & = &
\left. \frac{\nabla^N}{dw}\right|_{w=0}
\left. \frac{d}{dt}\right|_{t=s}g(t)g(s)^{-1}\gamma (w)\\ & = & \left.
\frac{\nabla^N}{dw}\right|_{w=0} X(\gamma (w))\\
 & =& \nabla^N X(\dot{\gamma}(0)).
\end{eqnarray*}

  By the last Lemma     $X$ is a normal vector field of 
$\exp^N (\bot_{c(s)}M)$. As $\exp^N (\bot_{c(s)}M)$ is totally geodesic in $N$, we see that $\nabla^N X
(\dot{\gamma}(0))$ is normal to $\exp^N (\bot_{c(s)}M )$. This means that
$\frac{\nabla^N}{dt}|_{t=s}g(t) _{\ast}v= \nabla^N X(\dot{\gamma}(0))$ is perpendicular to
$\bot_{c(s)}M$. Hence, we get   
$ g(t)_{\ast}  v = (\|^t_0 c)v$. \hfill $q.e.d.$

\section{Focalizations of equifocal submanifolds}

  Let $(G,K)$ be a  Riemannian symmetric pair of compact type such that    $
N=G/K$ is  simply connected. Let $N$ be   endowed with the
 metric that comes from the  Killing form of $\mathfrak{g}$. Let $
\tilde {M}\subset N$  be a
complete, connected, equifocal submanifold,   $M$  a
 focal manifold of $\tilde {M}$ and $\tilde {B}_{\xi}$ the tube
constructed in section 4 (compare (5) and (6)) (for
$\xi \in \bot M,\ \exp^N(\xi )\in \tilde{M}$). We will show that   $\tilde{M}$ and $\tilde{B}_{\xi}$   are equal. The
fact that
$\tilde{M}$  has abelian and compact normal bundle  (which was proved by  Terng and Thorbergsson, compare \cite{TT},
Theorem 6.15) and the normal parallel transport of
$\tilde{M}$  therefore is given by isometries of the ambient space
(compare Corollary \ref{milzehn}), is an important step in the proof.  

  So, let $\tilde{M}$ be an embedded, complete, connected, equifocal submanifold in $N$, 
$\tilde{\xi }:\tilde {M}\rightarrow\bot\tilde {M}$ a    parallel  normal vector field of $
\tilde {M}$ and $M:=\exp^N(\tilde{\xi }(\tilde {M}))$. 
Then,  $M$ is a submanifold of $N$  (compare page 3 in   \cite{TT}).   

\begin{lemma}\label{nag}
Let   $q\in M$ and  $\eta \in \bot_q  M$ with $z:= \exp^N (\eta ) \in \tilde{M}$. Then,
 $\mbox{\rm Hol}_{\eta}(M)\subset
\tilde{M}$.
\end{lemma}

  {\it proof}\hspace{0.3cm}   As the normal space of $\tilde{M}$ in $z$ is abelian  we
get $T_q (\exp^N (\bot_z \tilde{M}))\subset \bot_q M$ (compare Lemma \ref{drei}). By the fact that $T_z
(\mbox{Hol}_{\eta}(M)) =\exp^{\bot}_{\ast |\eta}(T_q M,0)$ and Lemma
\ref{drei} it follows that $\bot_z \tilde{M} \subset
\bot_z (\mbox{Hol}_{\eta} (M))$. Hence, $T_z (\mbox{Hol}_{\eta}(M)) \subset T_z \tilde{M}$. The proof then follows by  the
Theorem of Frobenius on foliations, because $\tilde{M}$ induces a  foliation on
$N$ in such a way that the parallel manifolds are equifocal again and $M$ is a focal manifold of each of them.  \hfill
$q.e.d.$

   {\bf Remark} \hspace{0.3cm} $M$ is not necessarily simply connected. Let $f: U\rightarrow M$ be the universal
covering, $\tilde{p}\in U$ and $\xi\in\bot_{f(\tilde{p})} M=\bot_{\tilde {p}}U$. Moreover, let $B_{\xi}$  be the tube
around $U$  defined in (5) and $ D_{\xi}$ the tube around $M$ defined  at the end of section 4. We already saw  
(compare Corollary
\ref{kyo})  
$\exp^ N(B_{\xi})=\exp^N(D_{\xi})=\tilde {B}_{\xi}$.

\begin{lemma}\label{hunviervierzig}
Let $\xi \in \bot_p M$  such that $\exp^N (\xi ) \in \tilde{M}$. Then, $\tilde{B}_{\xi}=\tilde{M}$.
\end{lemma}

  {\it proof}\hspace{0.3cm}  Let  
$\sigma :\tilde{M}\rightarrow M ,\ z \rightarrow \exp^N (\tilde{\xi}(z))$ be the canonical projection. $\sigma$ has
constant rank. Let $\sigma^{-1}(\{ q\} )\subset
\tilde{M}$  be the fibre over $q$ and  
\[ \tilde{\sigma}_q:=(\exp^N_q|\bot_qM)^{-1} (\sigma^{-1}(\{q\}))\subset\bot_qM.\]
Let $\hat{\sigma}_p$ be the connected component of $\tilde{\sigma}_p$ with $\xi \in \hat{\sigma}_p$.  In order to show
$\tilde{B}_{\xi} =\tilde{M}$ it is sufficient to show
$\hat{\sigma}_p =\hat{G}_p \xi$. Because  by Lemma \ref{nag} we then get $\tilde{B}_{\xi}=\tilde{M}$.

  We first want to show $\hat{G}_p \eta \subset \hat{\sigma}_p$ for all $\eta \in \hat{\sigma}_p$. Let us
therefore fix $\rho \in \tilde{\sigma}_q $ and $\tilde{A}:=
\exp^N (\bot_{\exp^N (\rho ) }\tilde{M})$. By the fact that $\tilde{A}$ is flat and totally geodesic in $N$ we get
\[  \bot_{\rho}
(\tilde{\sigma}_q ) =T_q (\exp^N (\tilde{A})),\]
where $\bot_{\rho} (\tilde{\sigma}_q )$ is the normal space of $\tilde{\sigma}_q$ in $\bot_q M$.   Let  $\nu \in
\bot_{\rho}(\tilde{\sigma}_q )$.
As $\rho$ and $\nu$ are contained in an abelian subspace  of $T_q N$ we get $\langle R_q^N (x,y)\rho ,\nu
\rangle =\langle R_q^N (\rho ,\nu )x,y\rangle =0$ for all $x,y \in
\bot_q M$. Thus,
$\tilde{R}_q(x,y)\rho \in  T_{\rho} (\tilde{\sigma}_q ) $ (for the definition of $\tilde{R}_q$ compare (1)).
 This can be done for any $q$.

   On the other hand we already saw   $\exp^N ((\|^1_0 c)\eta  )\in \tilde{M}$ for all
 curves $c:[0,1]\rightarrow M , c(0)=p$ and $\eta \in \hat{\sigma}_p$.  And so we get 
\begin{equation} \tau_c (
\tilde{R}_{c(1)})(x,y)\eta \in  T_{\eta} (\hat{\sigma}_p )  
\end{equation}
   for any 
$x,y \in \bot_p M$ and $\eta \in \hat{\sigma}_p$, $\nu \in  \bot_{\eta}(\hat{\sigma}_p ) $ (for the definition of 
$\tau_c ( \tilde{R}_{c(1)})$ compare (2)).

 Terng and Thorbergsson
showed in \cite{TT} that $\hat{\sigma}_p$ is isoparametric in $\bot_p M$ and   induces a (singular)  foliation. The
leaves of this foliation are mapped by $\exp^N$ to the leaves of the equifocal foliation in $N$.  A leaf
$D$ nearby $\hat{\sigma}_p$ is mapped to a regular leaf in $N$ (this means to an 
equifocal submanifold  in $N$). 

  Let $\tilde{R}, \bar{R} \in \hat{R}_p$ and $x,y,\bar{x},\bar{y} \in \bot_p M$ be chosen arbitrarily (for the definition
of $\hat{R}_p$ compare (3)). Let
$W
$ be the {\it fundamental vector field} corresponding to $\tilde{R}(x,y)$ (this is to say $W(\eta ) =\left.
\frac{d}{dt}\right|_{t=0} \exp^{\hat{G}_p}(t\tilde{R}(x,y)) \eta $ for all $\eta \in \bot_p M$)  and
$\bar{W}$ the fundamental vector field corresponding to $\bar{R}(\bar{x},\bar{y})$ respectively.  By (10) we get
$[W,\bar{W}]_{\eta}
\subset T_{\eta}(\hat{\sigma}_p )$ for all
$\eta
\in \hat{\sigma}_p$, where $[,]$ is the Lie bracket for vector fields. But $[W,\bar{W}]$ is the fundamental vector field
corresponding to $-[ \tilde{R}(x,y),\bar{R}(\bar{x},\bar{y})]_{{\cal L}_p}$. Hence, $T_{\eta}(G_p \eta )  \subset
T_{\eta} (\hat{\sigma}_p ).$ By the Theorem of Frobenius we get $G_p\eta\subset
\hat{\sigma}_p$.
 We already  saw  $\phi^{\ast}_p\eta\subset\hat{\sigma}_p$   (compare Lemma \ref{nag}). And so we get 
\begin{equation}
\hat{G}_p \eta \subset \hat{\sigma}_p.
\end{equation}

  Now, we want to prove $\hat{\sigma}_p \subset \hat{G}_p \eta $ for $\eta \in \hat{\sigma}_p$.  Let us therefore fix 
a curve
$\tilde{c}:{\RR} \rightarrow \hat{\sigma}_p, \tilde{c}(0)=\eta $. Let  
$\tilde{\eta}  \in \bot_{\exp^N (\eta )}\tilde{M}$ such that $\exp^N (\tilde{\eta})=p$ and   $\tilde{\eta}(s)$ be the
parallel normal vector field of $\tilde{M}$ along $\exp^N (\tilde{c})$ with $\tilde{\eta}(0)= \tilde{\eta}$. As the normal
 spaces are tori there exists a curve $k:{\RR} \rightarrow G$ such that 
$k(t)_{\ast}\tilde{\eta }(0)=\tilde{\eta }(t)$ (compare Corollary \ref{milzehn}).  But
$\exp^N (\tilde{\eta}(t))=p$ for all $t$. Hence, we get   $
k(t)p=p$ for all  $t$, this means $k(t)\in K_p$, the isotropy group of  $
p$.
It follows
\begin{equation}
 \hat{\sigma}_p \subset K_p \eta \subset T_p N.
\end{equation}

  Let $W\subset SO(\bot_p M)$ be the connected Lie subgroup whose   
Lie algebra is generated by   $\{ \tilde{R}_p (x,y)|\ x,y \in \bot_p M\}$. Then,
 $W\subset \hat{G}_p$ by construction.   Moreover, $\bot_{\eta} (W\eta ) \cap \bot_p M \subset \bot_{\eta} (K_p \eta ) $
for all $\eta  \in \bot_p M$ as one can see as follows: Let   
$v\in
 \bot_{\eta }(W\eta ) \cap \bot_p M$. For all $x,y \in \bot_p M$  we get $ \langle R^N (\eta
,v)x,y\rangle =\langle R^N (x,y)\eta ,v\rangle =0$
  which implies $ R^N (\eta ,v)=0,$
  because $N$ is a symmetric space with nonnegative sectional curvature   
(compare Lemma \ref{som}).
 Hence, $\eta$ and $v$ are contained in an abelian subspace of $T_
pN$. Therefore,   $v\in \bot_{\eta}(K_p\eta )$.
  
   We get  $ T_{\eta }(K_p
\eta ) \cap \bot_p M \subset T_{\eta }(W \eta )$.
 By $ T_{\eta }(W\eta )\subset T_{\eta }(\hat{G}_p \eta )\subset T_{\eta }(\hat{\sigma}_p )$
(compare (11))  we then get $ T_{\eta }(W
\eta )= T_{\eta }(\hat{G}_p
\eta )=T_{\eta }(\hat{\sigma}_p )$   (compare (12)). 

  Thus, $\tilde{B}_{\xi}$ and $\tilde{M}$ have the same dimension. As both manifolds are connected and complete  
we get  $\tilde{M}=\tilde{B}_{\xi}$.
\hfill $q.e.d.$

  By Lemma \ref{hundreidreissig} and the last Remark at the end of section 4.4 we can see that $
\hat {G}_p$ and the
normal parallel tranport of $M$    preserve the focal structure of $M$. 

\begin{lemma}\label{mar}
Let $(G,K)$ be a Riemannian symmetric pair of compact type, such that 
$N=G/K$ is simply connected. Let $N$ be endowed with the Killing metric. Let $\tilde{M} \subset N$ be an 
 equifocal submanifold    
  and $M $ a focal manifold. Then, $\hat{G}_p$ and the normal  parallel transport of $M$ preserve the focal structure of 
$M$. 
\end{lemma}

\begin{corollary}
We especially see that this is true for the singular orbits of hyperpolar actions.
\end{corollary}

  {\it proof}\hspace{0.3cm} It remains to show that the singular orbits are focal manifolds of the principal
orbits. This follows   by Corollary  \ref{milzehn}, because the elements in $
H$ of a given homogeneous equifocal
manifold
  $\tilde{M}=H/\tilde{H} \subset N$ induce the normal parallel translation of $\tilde{M}$. \hfill $q.e.d.$

  {\bf Remark}\hspace{0.3cm}  Let $M$ be a focal manifold of an equifocal one. The proof of Lemma
\ref{hunviervierzig} shows that  $\hat{G}_p \eta=W\eta$, where $W\subset SO(\bot _p M)$ is
the   Lie subgroup,
 whose Lie algebra is generated by  $\{ \tilde{R}_p
(x,y)|x,y \in \bot_p M\}$. We especially see $\phi_p^{\ast} \eta \subset W\eta$. This implies that the reduced normal
holonomy group $\phi^{\ast}_ p$ acts trivially on $V_0^p$. 
In Lemma
\ref{hunneunzehn} we saw $\phi^{\ast}_p|(V^p_1\oplus\ldots\oplus 
V_k^p)\subset G_p$. As $\phi^{\ast}_p$ acts   trivially on $V_0^p$ we get $
\phi^{\ast}_p\subset G_p$. And we get $\{ (\|^1_0 c) \hat{G}_p \xi |\ c:[0,1]\rightarrow M, c(0)=p \}= \{ (\|^1_0
c) W \xi |\ c:[0,1]\rightarrow M, c(0)=p \}=: D_{\xi}$.

\begin{theorem}\label{hunsechsvierzig}
Let $(G,K)$ be a Riemannian symmetric pair of compact type, such that 
$N=G/K$ is simply connected. Let $N$ be endowed with the Killing metric. Let $\tilde{M} \subset N$ be equifocal and $M$ a focal manifold. Let    $\xi\in\bot_pM$
with 
$\exp^N(\xi )\in\tilde {M}$. Then, we get $\tilde{B}_{\xi}= \exp^N (D_{\xi})=\tilde{M}$. Moreover, $W$ acts polarly
on $\bot_p M$.
\end{theorem}

  An interesting special case is the following:

\begin{theorem}\label{achtdreissig}
 Let $(G,K)$ be a Riemannian symmetric pair of compact type, such that 
$N=G/K$ is simply connected. Let $N$ be endowed with the Killing metric.
Let  $\tilde {M}\subset N$ be an equifocal submanifold which has some  parallel normal vector field   $
\tilde{\xi }:\tilde {M}\rightarrow\bot\tilde {M}$ such that  $\exp^
N(\tilde{\xi }(\tilde {M}))=\{p\}$, where $p=[K]$ without loss of generality. Then,
$\tilde{M}$ is an orbit of  $K$.
\end{theorem}

  {\bf Examples} \vspace{-0.5cm} \begin{enumerate}

\item {\it Inhomogeneous examples:}\hspace{0.3cm} Ferus, Karcher and M\"{u}nzner construct in \cite{FKM} series of
inhomogeneous, isoparametric families with codimension  1  in spheres. In order to do
this they use  the theory of  Clifford algebras and Clifford systems. Some of their examples in spheres of odd
dimension are invariant under the canonical $S^1$-action  on
$S^{2m+1}\subset\CC^{m+1}$ and  can be projected to inhomogeneous equifocal families in the 
complex projective space
 $\CC P^m$ (compare \cite{Wan}).  

\item {\it Homogeneous examples:} \hspace{0.3cm}The principal orbits of hyperpolar actions are equifocal (compare
\cite{TT}, Theorem 2.1). A very interesting class of hyperpolar actions was introduced by \cite{Her}. Let $
(G,K)$ and
$(G,\tilde{K})$ be two Riemannian symmetric pairs of compact type. Then, $K$  acts hyperpolarly on
$G/\tilde{K}$.
\begin{enumerate} 
\item Let for
example $G:=
\mbox{SU}(n+1),\ K:=S(U(n-k)\times U(k+1)) $ and $\tilde{K}:= S(U(n) \times U(1))$. Then, $S(U(n-k)\times U(k+1))$ acts
hyperpolarly on ${\CC}P^n =\mbox{SU}(n+1)/S(U(n)\times U(1))$. It is easy to check that ${\CC}P^k$ and ${\CC}P^{n-k-1}$
(both totally geodesically embedded into ${\CC}P^n$) are singular orbits of this action. Thus, the principal orbits are
{\it full tubes} around them. (Here we call $\exp^N (\bot^r (M))$ a {\it full tube} of $M \subset N$, where $\bot^r (M):=\{
v \in \bot M |\ \| v\| =r \}$.) This is so because ${\CC}P^n$ has rank 1.
\item Let $G$ and $\tilde{K}$ be  defined as before and $K:= SO(n+1)$. It is easy to check that
${\RR}P^n$ (totally geodesically embedded) and the complex quadric $Q^{n-1}$ are singular orbits. The principal orbits
again are full tubes around them. 

\end{enumerate}
  There exist hyperpolar actions (on symmetric spaces of rank 1) that are not of the
type described above. Takagi for example classified the cohomogeneity one actions on ${\CC}P^n$ (compare \cite{Tak}).
Among them one can find hyperpolar actions which are not of ``Hermann-type".
\end{enumerate}

  {\bf Remark}\hspace{0.3cm} The fact that full, compact,  irreducible isoparametric submanifolds of codimension
$\geq 3$ of the Euclidean space are orbits of $s$-representations  was firstly shown by  Thorbergsson
  (compare \cite{Tho}). Later  Olmos found another proof of the same fact (compare \cite{Ol2}). An important step in
his proof is the so-called ``Homogeneous  Slice Theorem". This Theorem says that 
the fibres of the projection
$\tilde {M}\rightarrow M$ of a compact, irreducible, isoparametric submanifold $
\tilde {M}\subset\RR^N$ onto a full
focal manifold   $M$ are homogeneous. They are orbits of some  polar action  (compare
\cite{HOT}). Recently, Heintze and Liu gave another proof of Thorbergssons Theorem - not only for
isoparametric submanifolds of the Euclidean space, but more generally
for isoparametric submanifolds of Hilbert spaces. They again used the
Homogeneous Slice Theorem (compare \cite{HL}).  Theorem \ref{hunsechsvierzig} of this paper is an analogue for equifocal
submanifolds of symmetric spaces of compact type. But our situation is more complicated than the situation of the
Homogeneous Slice Theorem, because the ambient space is more complicated.   Therefore we have to consider the curvature of
the ambient space in order to construct $D_{\xi}$.

  {\bf Open questions:}
\begin{enumerate}

\item Are equifocal submanifolds of codimension $\geq 2$ of irreducible symmetric spaces of compact type homogeneous? The
examples of Wang we considered before have codimension 1.

\item We already saw that there does not  exist more than one equifocal family to a given focal manifold. Is there a way
to conclude from  the geometry of the focal manifold whether the corresponding equifocal family is homogeneous or not? 

\item Which properties does the normal holonomy tube have?  In which cases do normal holonomy tube and the tube
$\tilde{B}_{\xi}$ constructed in this paper coincide?

\end{enumerate}

{\bf Acknowledgment}\hspace{0.3cm}  The results of this paper were obtained during my stay at the Universidad Nacional de
C\'ordoba in  Argentina with a  scholarship from DAAD (``DAAD-Doktoranden-Stipendium aus Mitteln des zweiten
Hoch\-schul\-son\-der\-pro\-gramms". I would like to thank DAAD for the financial support and the members of  the
Department of Ma\-the\-ma\-tics in C\'ordoba for their hospitality. In particular I would like to thank Carlos  Olmos for
the many fruitful discussions and valuable comments.

\vspace{3cm}

{\small Universit\"at zu K\"oln, Weyertal 86-90, D-50931 K\"oln}\\
{\small {\it E-mail address:} mbrueck@mi.uni-koeln.de}

\end{document}